\documentclass[12pt]{amsart}
\usepackage{amssymb}
\usepackage{amsmath}
\usepackage{amsfonts}
\usepackage[usenames]{color}
\usepackage{array}
\usepackage{color}

\makeatletter
 
 \@addtoreset{equation}{section}
\makeatother

\newtheorem{theorem}{Theorem}[section]
\newtheorem{lemma}[theorem]{Lemma}

\theoremstyle{definition}

\newtheorem{example}[theorem]{Example}

\newtheorem{corollary}[theorem]{Corollary}

\theoremstyle{remark}

\numberwithin{equation}{section}

\begin{document}

\author{Andrei Pavelescu}
\address{Department of Mathematics, Oklahoma State University, Stillwater, OK 74078-1030, USA}
\email{andrei.pavelescu@okstate.edu}

\date{\today}
\thanks{The author was partially supported by NSF grant
DMS 1001962.}

\title{Some Maximal Commutative Subrings of $M_n(D)$}
 
\maketitle

\begin{abstract}  We anwer a question of Weibel and
show that maximal commutative subgrings
of Artinian semisimple rings are direct products of local rings.
\end{abstract}

\section{Introduction}\vspace{.1in}
The motivation for the results in this paper lies with an open problem from Charles Weibel's online version of  "The K-book: An Introduction to Algebraic K-theory" \cite{W}.  In the second chapter of the book, Weibel generalizes the notion of $\tilde{K}_0(R)$ for non-commutative rings $R$, provided one knows that maximal commutative subrings of simple Artinian rings are finite products of 0-dimensional local rings.   We thank Weibel for asking the question. 

In this paper we prove a slightly stronger result. Let $D$ denote a division ring. Then:

\begin{theorem}  \label{main}  Let $R$ be a maximal commutative subring
of $M_n(D)$.   Then
$R$ is a direct product of at most $n$ 0-dimensional local rings and $J(R)=N(R)$.
\end{theorem}

As a direct consequence we get:

\begin{corollary}  If $R$ is a maximal commutative subring
of $M_n(D)$ containing no nilpotent elements, then $R$ is Artinian.
\end{corollary}

If $D$ is finite dimensional over $Z:=Z(D)$, then any subring
of $M_n(D)$ containing the center is Artinian.   If $[D:Z]$ is infinite dimensional,
we give examples of maximal commutative subrings $R$ with nilpotent elements 
which are and are not Artinian.

\section{Proof of Theorem \ref{main}} \vspace{.1in}

\begin{proof}
If $ n=1$ since for every element $x\in R$, $C_D(x)=C_D(x^{-1})$, $R$ is a field and the theorem holds.\\

Let $n\ge 2$.  If every non-zero element of $R\backslash N(R)$ is invertible, it follows that $R$ is a local ring and $N(R)=J(R)$ is the maximal ideal of $R$; thus $R$ is zero-dimensional. So we may assume that $R\backslash N(R)$ contains non-invertible elements.\\

Let $A\in R\backslash N(R)$ be a non-zero element of minimal rank $r$. By seeing $A$ as an element of $End_{D}(D^n)$, we conclude that $Ker(A)\cap Im(A)= \{0\}$ (Since otherwise the Rank+Nullity Theorem would imply $rank(A^2)<rank(A)$ and, as $A$ is not nilpotent, this would contradict the minimality of $r$).  Let  $\mathcal{B} = \{e_1,e_2,....,e_n\}$ be a basis for $D^n$ as a left vector space over $D$ with $e_i\in Im(A)$ for $1\le i \le r$ and $e_i\in Ker(A)$ if $i>r$. With respect to this basis,
\[A= \left( \begin{array}{ccc}
A'& \vline & 0 \\ \hline
0& \vline & 0
\end{array}\right ),\]\\ for some $r\times r$ invertible matrix $A'$. But then\\ \[R\subseteq C_{M_n(D)}(A)\subseteq\left ( \begin{array}{ccc}
M_r(D)& \vline & 0 \\ \hline
0& \vline & M_{n-r}(D)
\end{array}\right). \]\\ This implies that \\ \[R=\left ( \begin{array}{ccc}
R_1 & \vline & 0 \\ \hline
0& \vline & R_2 
\end{array}\right)\simeq R_1\times R_2.\]\\ Since $R$ is  maximal commutative, it follows that 
$R_1$ and $R_2$ are maximal commutative subrings of smaller matrix rings
over $D$. 
By induction we have $$R\simeq R_1\times R_2 \simeq \prod_1^kR_i,$$ with $R_i$ local zero dimensional rings. Also, $$N(R)\simeq N(R_1)\times N(R_2)\simeq J(R_1)\times J(R_2)\simeq J(R).$$

\end{proof}
Note that from the proof we get $N(R)^n=0$ and that $R$ is a direct
product of at most $n$ local rings. \\

\section{Examples}\vspace{.2in}

Let $D$ be a division ring and $n > 1$ a positive integer.  Set $S=M_n(D)$.  Let $E_{i,j}\in S$ denote the elementary matrix with 1 as the $(i,j)$ entry and zero everywhere else.  Let   $N \in S$ be the nilpotent matrix
consisting of one Jordan block in standard Jordan form with 0 on the main diagonal.  The following results are 
straightforward matrix computations.
\begin{lemma} 
\label{1} 
$C_S(N) =D\cdot I_n \bigoplus D\cdot (E_{1,2}+E_{2,3}+...+E_{n-1,n})\bigoplus...\bigoplus D\cdot(E_{1,n-1}+E_{2,n})\bigoplus D\cdot E_{1,n} $  and if $L$ is a maximal subfield of  $D$,
then $C_S(LN)= L\cdot I_n \bigoplus L\cdot (E_{1,2}+E_{2,3}+...+E_{n-1,n})\bigoplus...\bigoplus L\cdot(E_{1,n-1}+E_{2,n})\bigoplus D\cdot E_{1,n}$.
\end{lemma}
  
Let $M = E_{1,n}$.   Let $L$ be a maximal subfield of $D$. \\

\begin{lemma}
\label{2}  $C_S(DM)$ consists of all matrices with $(i,j)$ entry zero for $j=1, i > 1$ and
$i=n, j < n$ and with the $(1,1)$ and $(n,n)$ entries equal and in $Z$.  
\end{lemma}
Assume $[D:Z]= \infty$.
\begin{example}   Take  $R$ to be the centralizer of  $LN \cup \{M\}$,  then
$R$ consists of the banded upper triangular matrices with diagonal entries in $Z$,
the last entry arbitrary and all other entries in $L$.   By Lemmas \ref{1} and \ref{2}, $R$ is
maximal commutative.   However, $R$ is not Artinian (not even Noetherian) since any $Z$-subspace of $DE_{1,n}$
is an ideal of $R$ and since $D$ is infinite dimensional over $Z$. 
\end{example}

\begin{example}   Take $R$ to be the centralizer   of   $L  + LN$.  Then $R$ is contained in $C_S(L)
= M_n(L)$ and since $N$ is a regular element in the ring of matrices over a field, we have 
$ R = L[N]$ is maximal commutative.  Since $R$ is a finite dimensional $L$-algebra,
$R$ is Artinian (and contains nilpotent elements).
\end{example}


\begin{thebibliography}{999}

\bibitem{W}  C.  Weibel, The K-book: An Introduction to Algebraic K-theory, 
http://www.math.rutgers.edu/~weibel/Kbook.html

\end{thebibliography}
\end{document}